\documentclass[letterpaper, 10 pt, journal, twoside]{ieeetran}
\usepackage[utf8]{inputenc}
\usepackage{algcompatible}
\usepackage{algorithm}
\usepackage{cite}
\usepackage[caption=false,font=footnotesize]{subfig}
\usepackage{caption}
\usepackage{color}

\IEEEoverridecommandlockouts

\usepackage{graphicx} 
\usepackage{amsmath} 
\usepackage{amssymb} 
\usepackage{bm}

\newtheorem{theorem}{Theorem}
\newtheorem{lemma}{Lemma}
\newtheorem{definition}{Definition}
\newtheorem{corollary}{Corollary}

\newtheorem{assumption}{Assumption}
  
\DeclareMathOperator*{\argmin}{arg\,min}

\DeclareMathOperator{\proj}{proj}

\newcommand{\RR}{\mathbb{R}}
\newcommand{\NN}{\mathbb{N}}

\newcommand{\Reg}{\texttt{R}}
\newcommand{\XX}{\mathcal{X}}
\newcommand{\xx}{\mathbf{x}}
\newcommand{\xst}{\mathbf{x}^\star_t}

\newcommand{\xit}{\mathbf{x}_{i,t}}
\newcommand{\xjt}{\mathbf{x}_{j,t}}
\newcommand{\xsit}{\mathbf{x}_{i,t}^\star}
\newcommand{\xsitt}{\mathbf{x}_{i,t+1}^\star}
\newcommand{\xitt}{\mathbf{x}_{i,t+1}}
\newcommand{\yy}{\mathbf{y}}

\newcommand{\yit}{\mathbf{y}_{i,t}}
\newcommand{\yjt}{\mathbf{y}_{j,t}}
\newcommand{\yits}{\mathbf{y}_{i,t+1}^\star}

\newcommand{\pt}{\bm{\phi}_t}
\newcommand{\oyit}{\overline{\mathbf{y}}_{t}}

\newcommand{\yitt}{\mathbf{y}_{i,t+1}}
\newcommand{\Pij}{\mathbf{P}_{ij}}

\newcommand{\git}{\mathbf{g}_{i,t}}
\newcommand{\ogt}{\overline{\mathbf{g}}_{t-1}}

\newcommand{\vv}{\mathbf{v}}
\newcommand{\ww}{\mathbf{w}}
\newcommand{\ppi}{\bm{\pi}}

\newcommand{\at}{\mathbf{a}_t}
\newcommand{\ait}{a_{i,t}}
\newcommand{\sit}{s_{i,t}}
\newcommand{\nuit}{\nu_{i,t}}

\title{Dynamic and Distributed Online Convex Optimization for Demand Response\\ of Commercial Buildings}

\author{Antoine Lesage-Landry and Duncan S. Callaway
\thanks{This work was funded by the Natural Sciences and Engineering Research Council of Canada.}
\thanks{A. Lesage-Landry and D.S. Callaway are with the Energy \& Resources Group, University of California, Berkeley, CA, USA \{\texttt{alesagelandry,dcal}\}\texttt{@berkeley.edu}.}%
}

\begin{document}

\maketitle
\thispagestyle{empty}
\pagestyle{empty}

\begin{abstract}
We extend the regret analysis of the online distributed weighted dual averaging (\texttt{DWDA}) algorithm~\cite{hosseini2016online} to the dynamic setting and provide the tightest dynamic regret bound known to date with respect to the time horizon for a distributed online convex optimization (OCO) algorithm.
Our bound is linear in the cumulative difference between consecutive optima and does not depend explicitly on the time horizon. We use dynamic-online \texttt{DWDA} (\texttt{D-ODWDA}) and formulate a performance-guaranteed distributed online demand response approach for heating, ventilation, and air-conditioning (HVAC) systems of commercial buildings. We show the performance of our approach for fast timescale demand response in numerical simulations and obtain demand response decisions that closely reproduce the centralized optimal ones. 
\end{abstract}


\IEEEpeerreviewmaketitle

\section{Introduction}

\IEEEPARstart{D}{emand} response (DR) can provide an important part of the additional flexibility required to operate electric power systems with high penetration of renewables~\cite{callaway2011achieving,taylor2016power}.
Commercial and industrial buildings are an important class of thermostatically controlled loads offering flexibility that can be leveraged in DR~\cite{hao2014ancillary}. A building's heating, ventilation, and air-conditioning (HVAC) unit and specifically the air handler fan speed, can be temporarily altered to provide DR services on fast time-scales, e.g., frequency regulation~\cite{maasoumy2013flexibility,hao2014ancillary,lin2015experimental,kim2016technologies}. These services are required to ensure the stability and resiliency of modern grids~\cite{callaway2009tapping,mathieu2012stateTPS}.
 
In this work, we propose a distributed online convex optimization (OCO) approach for DR of commercial buildings. We extend the static regret analysis of the online distributed weighted dual averaging (\texttt{DWDA}) algorithm from~\cite{hosseini2016online} to the dynamic setting. We propose the dynamic-online \texttt{DWDA} (\texttt{D-ODWDA}) which provides an adequate performance guarantee for real-time DR. 
The \texttt{D-ODWDA} dynamic regret bound outperforms all previous distributed OCO algorithms and compares to non-distributed ones. 

Using a distributed OCO-based approach, we design a highly scalable and uncertainty-adaptive online approach for DR of commercial buildings. Buildings do not need to share their decision variables. Only weighted averages of the dual variables are transmitted, thus promoting privacy. Communication requirements and delays are reduced due to strictly local information exchange. We now review the literature related to our work. 

\paragraph*{Related work}

Several distributed OCO algorithms have been first designed~\cite{raginsky2011decentralized,tsianos2012distributed,yan2012distributed,mateos2014distributed,lee2016distributed,akbari2015distributed,koppel2015saddle,li2018distributed} based on a static regret analysis.
Our work focuses on the dynamic and distributed setting because it guarantees adequate performance for multi-period DR. There is a limited but important body of work on dynamic and distributed OCO. These include work on mirror updates~\cite{shahrampour2017online,shahrampour2017distributed}, adaptive search directions~\cite{nazari2019adaptive}, gradient-free methods~\cite{pang2019randomized} and time-varying constraints~\cite{sharma2020distributed}. This paper extends this body of work by using a distributed weighted dual averaging update~\cite{hosseini2016online,duchi2011dual} in the dynamic setting.  By doing so, we identify a stronger dynamic regret bound than prior studies.

Reference~\cite{wang2019development} surveyed HVAC-based DR approaches for commercial buildings. References~\cite{hao2014ancillary,lin2015experimental} proposed a controller for air handler fans to provide regulation. 
A hierarchical controller was presented in~\cite{vrettos2016experimental} using robust optimization and model predictive control.
A model predictive approach was used in~\cite{maasoumy2013flexibility,maasoumy2014model} to control the fan speed and the cooling and heating units. In~\cite{beil2016frequency}, the authors utilized the building's thermostat setpoint to control the fan speed and adjust the power consumption. Reference~\cite{hughes2016identification} formulated a virtual battery model for commercial buildings. In our work, we use distributed OCO for real-time DR of large aggregations of commercial buildings.

In this work, we make the following specific contributions:
\begin{itemize}
	\item We extend the regret analysis of~\cite{hosseini2016online} to the dynamic setting and prove, to the best of our knowledge, the tightest dynamic regret bound with respect to $T$ for distributed OCO algorithms. We present dynamic regret bounds and discuss the implication of our tighter bound.
	\item We propose a distributed online approach for commercial buildings equipped with HVAC and variable frequency drive-operated air handler fans for real-time DR. Our approach is highly scalable, promotes privacy and requires only local communications between the buildings. 
	We use computationally efficient decision-making updates thus tailoring our approach to real-time DR like frequency regulation.
\end{itemize}

\section{Preliminaries}
\label{ssec:prelim}
In this section, we present our notation and introduce the OCO framework formally.

\paragraph*{Notation}
We consider a time horizon $T \in \NN$ discretized into rounds $t$. We consider $n$ agents denoted by the subscript $i$. At round $t$, an agent makes a decision $\xit \in \XX$ with the objective of making the centralized, optimal decision, where $\XX \subseteq \RR^m$ is the decision set and $m\in \NN$.  
Let $f_t:\XX \mapsto \RR$ be the network loss function. Let $f_{i,t}:\XX \mapsto \RR$, a convex function, denote the local loss function of agent $i$. It is assumed to be differentiable for now, but this requirement is relaxed by Corollary~\ref{cor:subgradient}. The loss functions $f_t$ and $f_{i,t}$ are related by $f_t\left(\xx_t \right) = \frac{1}{n}\sum_{i=1}^n f_{i,t}\left(\xx_t \right)$. Finally, we consider the problem $\min_{\xx \in \XX} f_t(\xx)$, for $t=1,2,\ldots, T$, to be solved in an online and distributed fashion. In this setting, $f_t$ is only observed after round $t$.

We consider an undirected graph $\mathcal{G}=\left( \mathcal{V}, \mathcal{E} \right)$. The set of vertices $\mathcal{V} = \left\{1,2,\ldots,n \right\}$ represents the $n$ agents in the network. The edge set $\mathcal{E} \subseteq \mathcal{V} \times \mathcal{V}$ represents the undirected communication link between two agents. We define $\mathbf{P} \in \RR^{n^2}$ as the network matrix where $\Pij > 0$ if and only if there exists a communication channel between agent $i$ and $j$.

Let $\psi:\XX \mapsto \RR$ be a proximal function. We assume $\psi$ is 1-strongly convex, non-negative, and $\psi(\mathbf{0})=0$~\cite{duchi2011dual}. Based on this function, we have the following definition.
\begin{definition}[{\cite{duchi2011dual}}]
The $\psi$-regularized projection onto $\XX$ of a vector $\yy \in \RR^m$ with stepsize $\alpha$ is:
$\proj_\XX^{\psi}\left( \yy , \alpha \right) = \argmin_{\xx \in \XX}\left( \left\langle \yy, \xx \right\rangle + \frac{1}{\alpha}\psi(\xx) \right)$.
\label{def:proj}
\end{definition}
Let $V_T = \sum_{t=1}^T \left\| \xx_{t+1}^\star - \xx_t^\star \right\|$. The term $V_T$ characterizes the cumulative difference between consecutive optima. 

Lastly, we let $\langle \cdot,\cdot \rangle: \RR^m \times \RR^m \mapsto \RR$ be a scalar product. Let $\| \cdot \|$ denote a norm and $\| \cdot \|_\ast$ its associated dual norm. The dual norm is defined as: $\left\| \xx \right\|_\ast = \sup\left\{\left. \left\langle \xx, \yy \right\rangle \right| \yy \in \RR^m, \left\|\yy \right\| \leq 1 \right\}$ for $\xx \in \RR^m$.

\paragraph*{Background}

We now list our assumptions and discuss them briefly.
\begin{assumption}
\label{ass:strongly_connected}
The graph $\mathcal{G}$ is strongly-connected.
\end{assumption}

In other words, we assume that there exists a communication path linking all agents.
\begin{assumption}
All agents are connected to a minimum of two other agents and the network matrix $\mathbf{P}$ is row-stochastic, irreducible, and ergodic.
\label{ass:p}
\end{assumption}
It follows from Assumption~\ref{ass:p} that there exists a steady-state distribution $\bm{\pi} \in \left]0,1 \right[^n$ such that $\bm{\pi} = \bm{\pi} \mathbf{P}$ and $\mathbf{1}^\top \bm{\pi} = 1$ where $\mathbf{1}$ is the $n$-dimensional one vector.
\begin{assumption}
The set $\XX$ is compact and convex.
\label{ass:conv_comp}
\end{assumption}
Convexity and compactness of the decision set are standard assumptions in OCO~\cite{shalev2012online,hazan2016introduction}. 
Because the loss function is convex and its domain is compact, the loss function $f_{i,t}\left(\xx\right)$ is $L$-Lipschitz with respect to a norm $\|\cdot \|$ for all $i=1,2,\ldots, n$, $t=1,2,\ldots, T$ and $\xx \in \XX$~\cite{urruty1996convex}. Consequently, $\left\|\nabla f_{i,t}\left(\xx\right)\right\|_\ast \leq L$ for all $i$, $t$, and $\xx \in \XX$~\cite[Lemma 2.6]{shalev2012online}. Lastly, the problem is assumed to be feasible at all $t$ with the optimum denoted by $\xsit$, $i=1,2,\ldots,n$ and all $t$. From Definition~\ref{def:proj}, there exists $Y < \infty$ such that $Y \geq \left\| \yy \right\|$ for all $\yy \in \mathcal{Y}^\star_{i,t}=\left\{\left. \yy \in \RR^m \right| \xsit = \proj_\XX^{\psi}\left( \yy , \alpha \right) \right\}$, 
$i=1,2,\ldots, n$, and $t=1,2,\ldots,T$.

\begin{assumption}
The time horizon $T \in \NN$ is known.
\label{ass:time_horizon}
\end{assumption}
This is a standard assumption in OCO~\cite{shalev2012online,hazan2016introduction,koppel2015saddle,shahrampour2017distributed}. If the time horizon $T$ is unknown a priori, i.e., the forecaster does not know when the control process will end, the forecaster can initialize the problem with time horizon $T'\in \mathbb{N}$. If $t=T$ is not reached within $T'$, the algorithm can be reinitialized for the next $T'$ rounds, and this process can be repeated until $t=T$ is reached.

The performance of an OCO algorithm is characterized by its regret. In this work, we evaluate the performance using the dynamic regret: $\Reg_T = \sum_{t=1}^T \left( f_t\left(\xx_t\right) - f_t\left(\xst\right) \right)$, where $\xst \in \argmin_{\xx \in \XX} f_t\left(\xx\right)$, the round optimum. The dynamic regret differs from the static regret in which the best-fixed decision in hindsight is used as comparator: $\xst = \xx^\star \in \argmin_{\xx\in\XX} \sum_{t=1}^T f_t\left(\xx\right)$ for all $t=1,2,\ldots,T$ instead of the round optimum. The static regret offers a sufficient performance guarantee for certain applications like the estimation of a static vector via sensor networks~\cite{akbari2015distributed,hosseini2016online} or the training of support vector machines for security breaches~\cite{koppel2015saddle}. It is an inadequate performance metric for dynamic or time-varying problems, e.g., localization of moving targets~\cite{shahrampour2017distributed,lesagelandry2019second}, or setpoint tracking like real-time DR~\cite{lesage2018predictive}. For these problems, a dynamic regret analysis is required to offer a suitable performance guarantee.

In distributed OCO, several regret definitions are used. The network or coordinated dynamic regret~\cite{koppel2015saddle,lee2016distributed} compares the agents' decisions to the centralized round optimum, $\xst \in \argmin_{\xx \in \XX} \frac{1}{n}\sum_{i=1}^n f_{i,t} \left( \xx \right)$. The network dynamic regret is:
\begin{align}
\Reg_T &= \sum_{t=1}^T \left( \frac{1}{n}\sum_{i=1}^n f_{i,t}\left(\xit\right) - f_t \left( \xst \right)\right)\label{eq:recall_regret}
\end{align}
An algorithm with a sublinear dynamic regret bound will perform, at least on average, as well as the centralized round optimum as the time horizon increases~\cite{shalev2012online,hazan2016introduction}. 
The local dynamic regret for agent $j$ compares the decision played by an agent as if it was implemented by all to the centralized optimum~\cite{mateos2014distributed,koppel2015saddle,lee2016distributed,hosseini2016online,mokhtari2016online}. It is given by:
\begin{align}
\Reg_T(j) &= \sum_{t=1}^T \left( f_t\left(\xx_{j,t}\right) - f_t \left( \xst \right)\right) \label{eq:agentj_regret}
\end{align}
If the local regret is sublinearly bounded, then all agents will play, on average, as well as a centralized round optimum as $T$ increases. Thus, the agents have learned from their neighbours~\cite{koppel2015saddle,lee2016distributed}.

\section{Dynamic-\texttt{ODWDA}}
\label{sec:dodwda}
The \texttt{D-ODWDA} algorithm is presented in Algorithm~\ref{alg:d-odwda}. Our work is based on earlier \texttt{DWDA} method~\cite{hosseini2016online} but focuses on dynamic rather than static regret.  This is a more general metric and the proof strategy differs from earlier work. The \texttt{D-ODWDA} update is the following~\cite{hosseini2016online,duchi2011dual}:
\begin{align}
\yitt &= \sum_{j=1}^n \Pij \yy_{j,t} + \nabla f_{i,t} \left( \xit \right) \label{eq:update_y}\\
\xitt &= \proj_\XX^{\psi}\left( \yitt , \alpha \right) \label{eq:update_x}
\end{align}
The update is similar to~\cite{hosseini2016online} except for the constant step size $\alpha$. This parameter enables us to bound the difference between asynchronous regularized projections and guarantee dynamic performance. Let $\oyit,\pt, \ogt \in \RR^m$. These vectors are respectively the weighted average of the dual variables, $\oyit = \sum_{i=1}^n \pi_i \yit$, its regularized projection, $\pt = \proj_\XX^\psi \left(\oyit, \alpha_t \right)$, and the weighted average of the gradients, $\overline{\mathbf{g}}_t = \sum_{i=1}^n \pi_i \nabla f_{i,t} \left( \xit \right)$.

\setlength{\textfloatsep}{0pt}
\begin{algorithm}[h]
\begin{algorithmic}[1]
\STATEx \textbf{Parameters:} $T$, $\beta$, $\psi\left( \xx \right)$, $\mathbf{P}$. 
\STATE  Set $\mathbf{y}_{i,0} = \mathbf{0}$, $\xx_{i,0} \in \XX$ for all $i$, and $\alpha =\frac{\beta}{T}$

\FOR{$t = 0, 1,2, \ldots, T$}
\STATE Implement $\xx_{i,t}$ and observe the loss functions $f_{i,t}(\xx_{i,t})$
\STATE For all $i=1,2,\ldots, n$, update  $\yitt$ and $\xitt$:
\begin{align*}
\yitt &= \sum_{j=1}^n \Pij \yy_{j,t} + \nabla f_{i,t} \left( \xit \right) \\
\xitt &= \proj_\XX^{\psi}\left( \yitt , \alpha \right)
\end{align*}
\ENDFOR
\end{algorithmic}
\caption{Dynamic-\texttt{ODWDA}}
\label{alg:d-odwda}
\end{algorithm}

\subsection{Technical lemmas}

We present two lemmas which are then used in the regret analysis. Their proofs are provided in Appendices~\ref{app:error_y} and~\ref{app:yit}.

\begin{lemma}[{\cite[Lemma 2]{duchi2011dual}}]
\label{lem:proj_ineq}
Let $\vv, \ww \in \RR^d$, where $d \in \NN$ and $\eta > 0$, then
$
\left\|\proj_\XX^{\psi}\left( \vv , \eta \right) - \proj_\XX^{\psi}\left( \ww , \eta \right)\right\| \leq \eta \left\|\vv - \ww \right\|_\ast.
$ 
\end{lemma}

\begin{lemma}
\label{lem:dual_avg_err}
The following bound holds: $\left\| \oyit - \yit \right\|_{\ast} \leq \frac{nL}{\gamma\left(1-\gamma^{\frac{1}{\nu}}\right)}  + 2L$,
for all $i$ and $t$ and where $\gamma \in [0,1[$ and $\nu \in \mathbb{N}$ such that $\nu \geq 1$ are defined as in~\cite[Theorem 1]{anthonisse1977exponential}.
\end{lemma}

\begin{lemma}{}
\label{lem:norm_yit}
Let $p=\max_{i,j} \Pij$. The dual norm of $\oyit$ is bounded above and $
\left\| \oyit\right\|_\ast \leq \frac{n^2 L}{1-p}  + L,$ for all $t=1,2,\ldots, T$.
\end{lemma}

\subsection{Regret analysis}
\label{sub:regret_analysis}

We now present our main results. We show that \texttt{D-ODWDA} has a sublinear network dynamic regret bound. We then prove that a similar bound holds for the local dynamic regret and for sub-gradient-based updates.

\begin{theorem}[Network dynamic regret bound]
\label{thm:dyn_regret}
Suppose Assumptions~\ref{ass:strongly_connected}-\ref{ass:time_horizon} are met and let $\mathbf{y}_{i,0}= \mathbf{0}$ and $\alpha_t = \frac{\beta}{T}$ for all $t=1,2,\ldots, T$ and $i=1,2,\ldots,n$. The network dynamic regret of \texttt{D-ODWDA} is then bounded above by:
\begin{align*}
\Reg_T &\leq \beta L^2 \left(\frac{n}{\gamma\left(1-\gamma^{\frac{1}{\nu}}\right)}  + 2 \right)\\
&\quad + \beta L \left(\frac{n^2 L}{1-p} + L + Y\right) + L V_T.
\end{align*}
\end{theorem}

\begin{IEEEproof}
The objective function $f_{i,t}$ is convex for all $t$ and $i$, thus for $\xit, \xsit \in \RR^m$, the following inequality holds:
\begin{equation}
f_{i,t}\left( \xit \right) -f_{i,t}\left(\xsit\right) \leq  \left\langle \nabla f_{i,t} \left( \xit \right), \xit - \xsit  \right\rangle \label{eq:by_convexity}
\end{equation}
Using~\eqref{eq:by_convexity} we upper bound~\eqref{eq:recall_regret} and re-write the regret as
\begin{align}
\Reg_T &\leq \sum_{t=1}^T \frac{1}{n}\sum_{i=1}^n \left\langle \nabla f_{i,t}\left(\xit\right),\xit - \xsit \right\rangle \nonumber\\
&\leq \sum_{t=1}^T \frac{1}{n}\sum_{i=1}^n L \left\|\xit - \pt\right\| + L\left\|\pt - \xsit \right\|, \label{eq:inter_2sums}
\end{align}
where the last inequality follows from Lipschitz continuity of $f_{i,t}$. Using Lemma~\ref{lem:proj_ineq}, we upper bound the first term of~\eqref{eq:inter_2sums}:
\begin{equation}
\Reg_T\leq \sum_{t=1}^T \frac{1}{n}\sum_{i=1}^n L \alpha_t \left\|\oyit - \yit \right\|_\ast + L \left\|\pt - \xsit \right\|. \label{eq:after_ae}
\end{equation}
We re-express the last term of~\eqref{eq:after_ae} as:
\begin{align*}
L \left\|\pt - \xsit \right\| &= L \left\|\pt - \xsit + \xsitt - \xsitt\right\| \nonumber\\
&\leq L \left\|\pt - \xsitt \right\| + L\left\| \xsitt - \xsit \right\|
\end{align*}
where we used the triangle inequality. By assumption, $\alpha_t = \alpha$ for all $t=1,2,\ldots T$. Thus, Lemma~\ref{lem:proj_ineq} can be used on terms with different time indices. This leads to
\begin{equation}
L \left\|\pt - \xsit \right\| \leq L \alpha \left\|\oyit - \yits \right\| + L\left\| \xsitt - \xsit \right\| \label{eq:looking_for_vt}
\end{equation}
where $\yits \in \mathcal{Y}^\star_{i,t+1}$.
We now upper bound $\left\|\oyit - \yits \right\|$ of~\eqref{eq:looking_for_vt}. Using the triangle inequality and Lemma~\ref{lem:norm_yit}, we obtain $\left\|\oyit - \yits \right\|
\leq \frac{n^2 L}{1-p} + L + \left\|\yits \right\|$. Because $Y \geq \left\|\yits \right\|$ for all $i$ and $t$, and we have $\left\|\oyit - \yits \right\| \leq \frac{n^2 L}{1-p} +L + Y.$ We re-write~\eqref{eq:looking_for_vt} as
\begin{equation}
L \left\|\pt - \xsit \right\| \leq \alpha L \left(\frac{n^2 L}{1-p} + L + Y\right) + L\left\| \xsitt - \xsit \right\|.\label{eq:presque}
\end{equation}
We invoke Lemma~\ref{lem:dual_avg_err} and use~\eqref{eq:presque} to upper bound both terms of~\eqref{eq:after_ae}. By definition, $\xsit = \xst$ for all $i=1,2,\ldots, n$ and recall that $V_T = \sum_{t=1}^T \left\| \xx^\star_{t+1} - \xst \right\|$. This leads to
\begin{align*}
\Reg_T&\leq \sum_{t=1}^T \alpha L^2 \left(\frac{n}{\gamma\left(1-\gamma^{\frac{1}{\nu}}\right)}  + 2 \right) \\
&\quad+ \sum_{t=1}^T \alpha L \left(\frac{n^2 L}{1-p} + L + Y\right) + L V_T.
\end{align*}
Setting $\alpha = \beta/T$ completes the proof.
\end{IEEEproof}

Consequently, $\Reg_T \propto O\left(1+ V_T\right)$, and thus the regret is sublinear if $V_T < O\left(T\right)$. We now extend Theorem~\ref{thm:dyn_regret} to the local dynamic regret.

\begin{corollary}[Local dynamic regret of agent $j$ bound]
Suppose Assumptions~\ref{ass:strongly_connected}-\ref{ass:time_horizon} hold and let $\mathbf{y}_{i,0} = \mathbf{0}$ and $\alpha_t = \frac{\beta}{T}$ for all $t=1,2,\ldots, T$ and $i=1,2,\ldots,n$, then the local dynamic regret of agent $j$ when using \texttt{D-ODWDA} is bounded above by Theorem~\ref{thm:dyn_regret}'s bound.
\label{cor:agent_regret}
\end{corollary}

\begin{IEEEproof}
We proceed similarly to Theorem~\ref{thm:dyn_regret}'s proof with $\xjt$ instead of $\xit$ in the local dynamic regret for agent $j$~\eqref{eq:agentj_regret}. This leads to:
\begin{align*}
\Reg_T(j) \leq \sum_{t=1}^T \frac{1}{n} \sum_{i=1}^n \big( L \alpha_t \left\| \oyit - \yjt \right\| + L \left\|\pt - \xsit  \right\|\big).
\end{align*}
We observe from Lemma~\ref{lem:dual_avg_err} that the first argument of the double sums can be bounded by a constant. The corollary then follows from Theorem~\ref{thm:dyn_regret}'s proof at~\eqref{eq:after_ae} and onwards.
\end{IEEEproof}
If $V_T < O(T)$, the local regret is sublinear. The time-averaged regret thus decreases as $T$ increases and the agent's decisions are similar to the centralized optima, on average.

The dynamic regret bounds presented in Theorem~\ref{thm:dyn_regret} and Corollary~\ref{cor:agent_regret} have a tighter dynamic regret bounds than all other distributed OCO algorithms we are aware of~\cite{shahrampour2017online,shahrampour2017distributed,pang2019randomized,nazari2019adaptive,sharma2020distributed}. Our $O \left( 1 + V_T \right)$ bounds have a smaller order of dependence on $T$, only in the $V_T$ term, than the algorithm which had previously the tightest proved bound~\cite{shahrampour2017distributed}. Reference~\cite{shahrampour2017distributed}'s bound is $O\left(\sqrt{1 + TV_T} \right)$ or $O\left(\sqrt{T} \left(1 + V_T \right) \right)$ with and without prior knowledge of $V_T$, respectively. Our improved results may be due, in part, to the dual weighted averaging update which was shown to achieve high performance in offline optimization~\cite{duchi2011dual} and to the slightly stronger assumption on the network (Assumption~\ref{ass:p}) where we assumed that each agent is connected to at least two agents in addition to the standard network assumptions~\cite{shahrampour2017distributed}. Unlike~\cite{shahrampour2017distributed} which requires $V_T < O\left( \sqrt{T} \right)$ to yield sublinear regret bounds, our results are sublinear for $V_T < O\left( T \right)$ and thus hold for a larger family of time-varying optimization problems. Finally, the $O\left( 1 + V_T \right)$ bound of \texttt{D-ODWDA} is of the same order as the tightest known dynamic regret bound for any non-distributed algorithm~\cite{mokhtari2016online}.

Let $\partial f_{i,t}\left(\xx\right)$ be the set of sub-gradients of $f_{i,t}$ at $\xx$. To conclude, we have the following corollary.

\begin{corollary}[Sub-gradient-based \texttt{D-ODWDA}]
Theorem~\ref{thm:dyn_regret} and Corollary~\ref{cor:agent_regret} hold
when sub-gradients are used instead of the loss functions' gradients.
\label{cor:subgradient}
\end{corollary}

\begin{IEEEproof}
Let $\git(\xx_t) \in \partial f_{i,t}\left(\xx_t\right)$ be a sub-gradient of $f_{i,t}$ at $\xx$. We remark that the same gradient bounds hold for $\left\|\git(\xx_t)\right\|$ and $\left\|\git(\xx_t)\right\|_\ast$~\cite[Lemma 2.6]{shalev2012online}. Next, observe that~\eqref{eq:by_convexity} holds when substituting $\git(\xx_t)$ to $\nabla f_{i,t} \left( \xx_t\right)$ by the definition of a sub-gradient. The regret bound follows from Theorem~\ref{thm:dyn_regret} where $\git(\xx_t)$ is used instead of $\nabla f_{i,t} \left( \xx_t\right)$.
\end{IEEEproof}

\section{Distributed online demand response}
\label{sec:dist_setpoint} 
We now formulate a distributed online demand response approach for commercial buildings based on \texttt{D-ODWDA}. The buildings modulate in real-time their air handler's speed~\cite{hao2014ancillary} to increase or decrease their electric power consumption and provide DR services. Specifically, we consider real-time power setpoint tracking with flexible loads. Solving the problem in a distributed fashion allows for our approach:
(i) to be highly scalable as each load computes their low-dimensional control, (ii) to reduce the communication requirement and concurrently, to minimize unreliable or corrupted communication issues between a centralized decision-maker and the buildings, and (iii) to promote privacy as power adjustment variables are never communicated. Only indirect information, $\yit$, about the loads is communicated to their neighbours.

\subsection{Formulation}
\label{ssec:formulation}
Our model is based on~\cite{lesage2018setpoint},~\cite{yang2016distributed}, and~\cite{hao2014ancillary} for, respectively, the setpoint tracking and the distributed optimization, and the commercial building DR settings. Each building can adjust the speed of their fan on a short timescale leading to an adjustment $\ait$ to their nominal power consumption. Reference~\cite{hao2014ancillary} showed that for a given setpoint bandwidth, the power consumption of the HVAC can be approximatively considered as linearly dependent with the fan speed adjustment. 
Let $\ait \in \left[\underline{a}_i , \overline{a}_i \right]$ and $\at = \begin{pmatrix} a_{1,t} & a_{2,t} & \ldots & a_{n,t} \end{pmatrix}^\top$ where $\underline{a}_i$ and $\overline{a}_i$ are, respectively, the minimum and maximum adjustment load $i$ can provide and $-\underline{a}_i$, $\overline{a}_i >0$. Let the decision set be $\mathcal{A} = \bigcup_{i=1}^n \left[\underline{a}_i , \overline{a}_i \right]$. The objective of the DR aggregator is to dispatch loads such that their total power consumption adjustment meets the setpoint $s_t \in \RR$ while minimizing the adjustment required from each building. For this purpose, we use the squared $\ell_2$-norm which will penalize large deviations from the building's nominal operation. The DR online optimization problem is:
\begin{equation}
\min_{\at \in \mathcal{A}} \left\| \at \right\|_2^2 \qquad \text{s.t.} \quad \sum_{i=1}^n \ait = s_t.
\label{eq:setpoint_tr}
\end{equation} 
Let $\nu_t \in \RR$ be the dual variable associated to the equality constraint of~\eqref{eq:setpoint_tr}. The dual problem of~\eqref{eq:setpoint_tr} is
\begin{equation}
\max_{\nu_t \in \RR} \sum_{i=1}^n \Gamma_i\left( \nu_t \right) - \nu_t s_t, \label{eq:dual_of}
\end{equation}
where $\Gamma_i\left(\nu_t \right) = \min_{\ait \in \left[\underline{a}_i , \overline{a}_i \right]} c_i \ait^2 + \nu_t \ait$. In their current form, neither~\eqref{eq:setpoint_tr} nor~\eqref{eq:dual_of} are distributed problems. We follow~\cite{yang2016distributed}'s approach to obtain an equivalent distributed online optimization problem. This approach relies on solving the dual problem using virtual setpoints~\cite{yang2016distributed}. Let $\sit \in \RR$ such that $\sum_{i=1}^n \sit = s_t$ for all $t$ be the virtual setpoints. We consider the associated distributed online dual problem:
\begin{equation}
\min_{\nuit \in \RR} - \Gamma_i\left(\nuit \right) + \nuit \sit,
\label{eq:dual}
\end{equation}
for all loads $i=1,2,\ldots, n$ and where the dual variable $\nu_t$ is decoupled into $n$ local dual variables $\nu_{i,t}$. The primal variable $\ait$ can then be computed in each round from $\nu_{i,t}$. We use \texttt{D-ODWDA} on~\eqref{eq:dual}. 
Each building computes their local adjustment to track the setpoint as follows. In each round $t$, building $i$ implements $a_{i,t} = \min \left\{ \max \left\{ -\frac{\nu_{i,t}}{2} ,\underline{a}_i\right\},\overline{a}_i \right\}$ and observes the outcome of the decision. The round concludes with the building updating $\nu_{i,t+1}$ using~\eqref{eq:update_y} and~\eqref{eq:update_x}.

By Corollary~\ref{thm:dyn_regret}, solving the online problem~\eqref{eq:dual} using \texttt{D-ODWDA} will lead to $\nu_{i,t} = \nu^\star_t$ for all $i$, at least on average as $T$ increases. By definition, $\nu^\star_t = \argmin_{\nu \in \RR} \sum_{i=1}^n \left(- \Gamma_i\left(\nu \right) + \nu \sit \right) \equiv \argmin_{\nu \in \RR} - \sum_{i=1}^n  \Gamma_i\left(\nu \right) + \nu s_t$, the dual optimum~\eqref{eq:dual_of}. 

By strong duality and strong convexity, it follows that as the time horizon increases, $\ait = a_{i,t}^\star$ at least on average as well and the buildings will implement the optimal adjustment dispatch on average. We note that because of strong duality and primal feasibility, there exists a convex and compact set $\mathcal{X} \subset \RR$ such that $\nuit \in \mathcal{X}$ for all $i$ and $t$ and all assumptions of Corollary~\ref{thm:dyn_regret} are met. We do not compare our approach to other OCO algorithms in this section. Regret analysis results are (i) sufficient conditions and (ii) do not characterize individual round performance but rather bound worst-case performance. Thus, a comparison would neither confirm nor inform our results.

\subsection{Numerical example}
\label{ssec:numerical}
We consider $4$-second frequency regulation rounds and a time horizon $T=1000$ equivalent to $66.66$ minutes. We let $n=5$ to better visualize the behavior of the distributed algorithm. For loads $i=\{3, 4, 5 \}$, we sample the maximum and minimum power adjustment capacity, $\overline{a}_i$ and $\underline{a}_i$, uniformly in $[2,3]$ and $[-3,-2]$ kW. We set the capacity to be between $\pm 0.5$kW and $\pm 0.75$kW for buildings $1$ and $2$. We assume that each agent is connected to their neighbours, e.g. load $1$ to $2$, $2$ to $3$ and $5$ to $1$, thus meeting Assumption~\ref{ass:p}. The setpoint to track is $s_t = s_{t-1} + \sigma\left( -1  \right)^{b_t}/\sqrt{t}$ where $b_t \sim \mathrm{Bernoulli}(0.5)$, $\sigma=2$ kW and $s_0 = 0$ kW. The virtual setpoints are set to $s_{i,t} = s_t / n$ for all $i$ and $t$. We set $\psi\left(\xx\right) = \frac{1}{2}\left\| \xx \right\|^2_2$ and let $\beta = 200$.

\begin{figure*}[!t]
\centering
\subfloat[Load dual variable vs. centralized dual optimum]{\includegraphics[width=0.66\columnwidth]{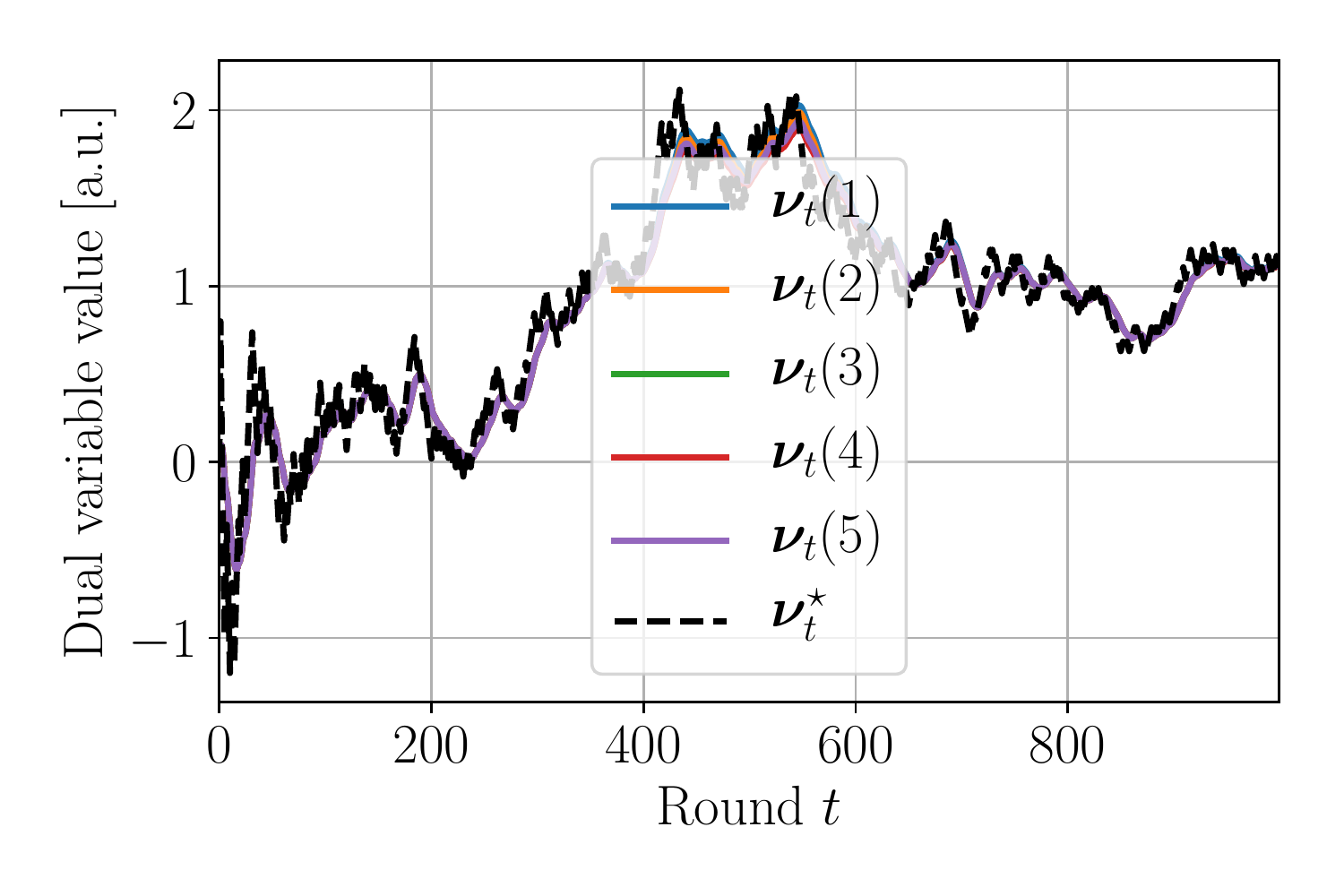}
\label{fig:agent_dual}}
\subfloat[Aggregated setpoint tracking]{\includegraphics[width=0.66\columnwidth]{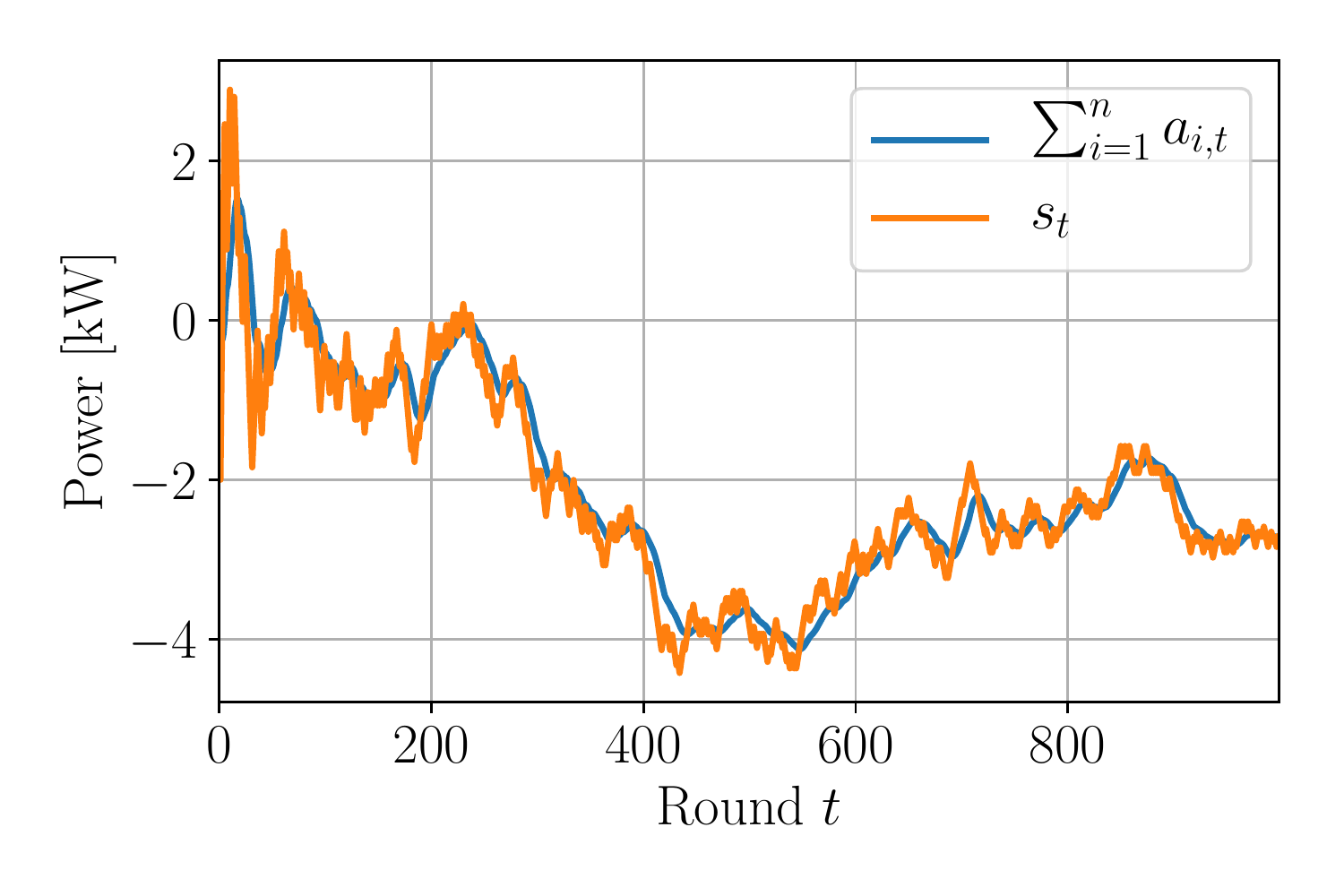} 
\label{fig:setpoint_agg}}
\subfloat[Local setpoint tracking]{\includegraphics[width=0.66\columnwidth]{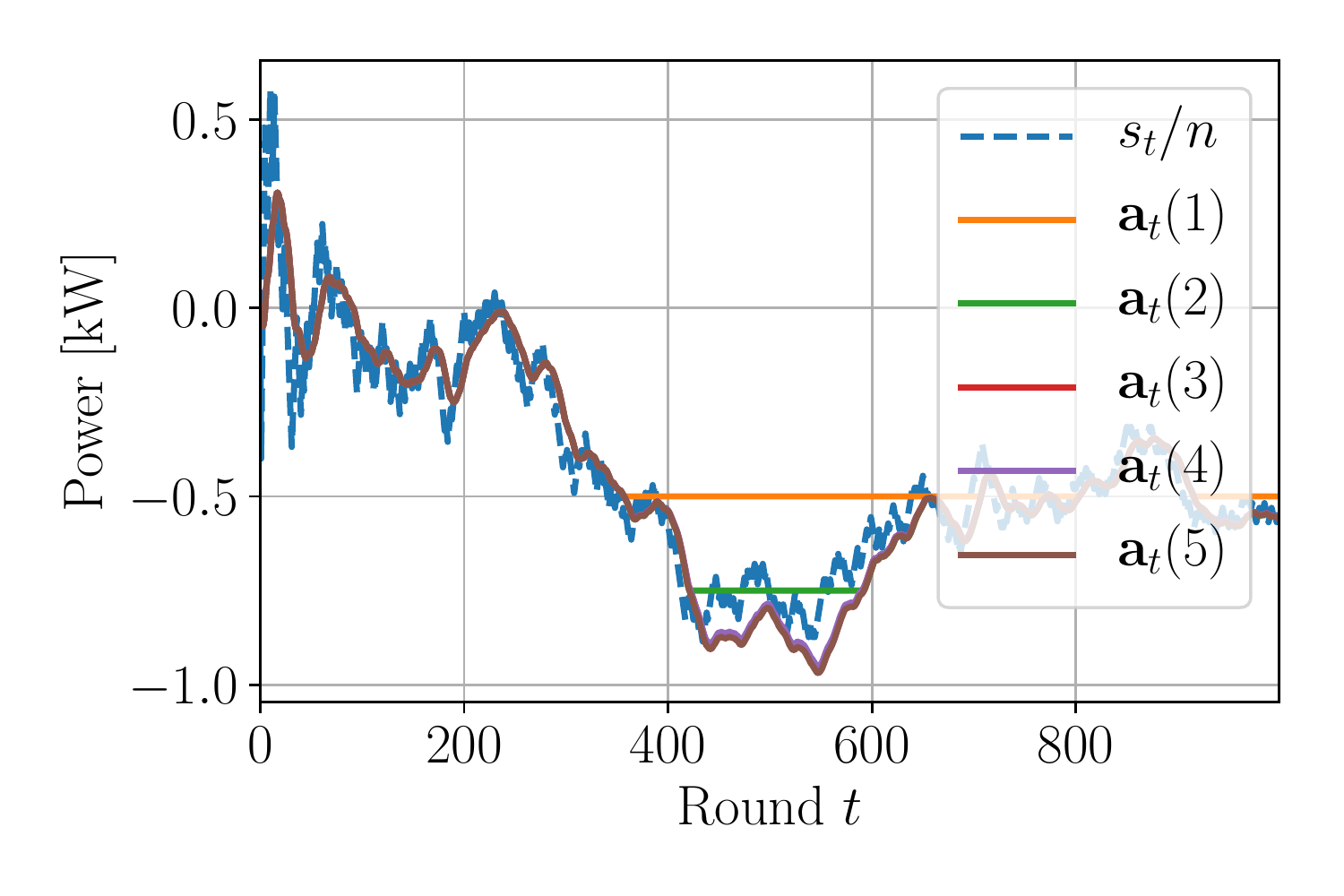}
\label{fig:agent_setpoint}}

\caption{Distributed setpoint tracking using \texttt{D-ODWDA}}
\label{fig:setpoint}
\vspace{-0.5cm}
\end{figure*}

Figure~\ref{fig:setpoint} presents the performance of our \texttt{D-ODWDA} for DR. Figure~\ref{fig:agent_dual} compares the load's dual variable $\nu_{i,t}$ to the centralized optimal value computed from~\eqref{eq:setpoint_tr} in hindsight. Figure~\ref{fig:agent_dual} shows that the dual variables computed by each building using \texttt{D-ODWDA} are similar to the centralized problem dual optimum, $\nu_t^\star$, with a relative difference at $T$ of $1.7\%$, $0.8\%$, $0.5\%$, $0.5\%$, $0.8\%$, for building $i=1,2,\ldots, 5$. The approach, therefore, computes power adjustments that are closely related to the centralized optimal adjustment. 
Figure~\ref{fig:setpoint_agg} presents the setpoint tracking performance from all buildings in the network and shows that our approach can adequately track the time-varying regulation setpoint. Figure~\ref{fig:agent_setpoint} presents local setpoint tracking and the virtual setpoint. We note that loads $1$ and $2$ cannot track their virtual setpoint at all rounds because of their limited adjustment capacity. During these instances, the other loads increase their contribution so that the setpoint $s_t$ is matched. This is shown in Figure~\ref{fig:agent_setpoint} in which loads $3$-$5$ have higher curtailment dispatched than their virtual setpoints. 
No centralized entity intervenes and only the exchange of $\yy_{i,t}$ variables suffice. 

\begin{figure}[tb]
  \centering
  \includegraphics[width=0.85\columnwidth]{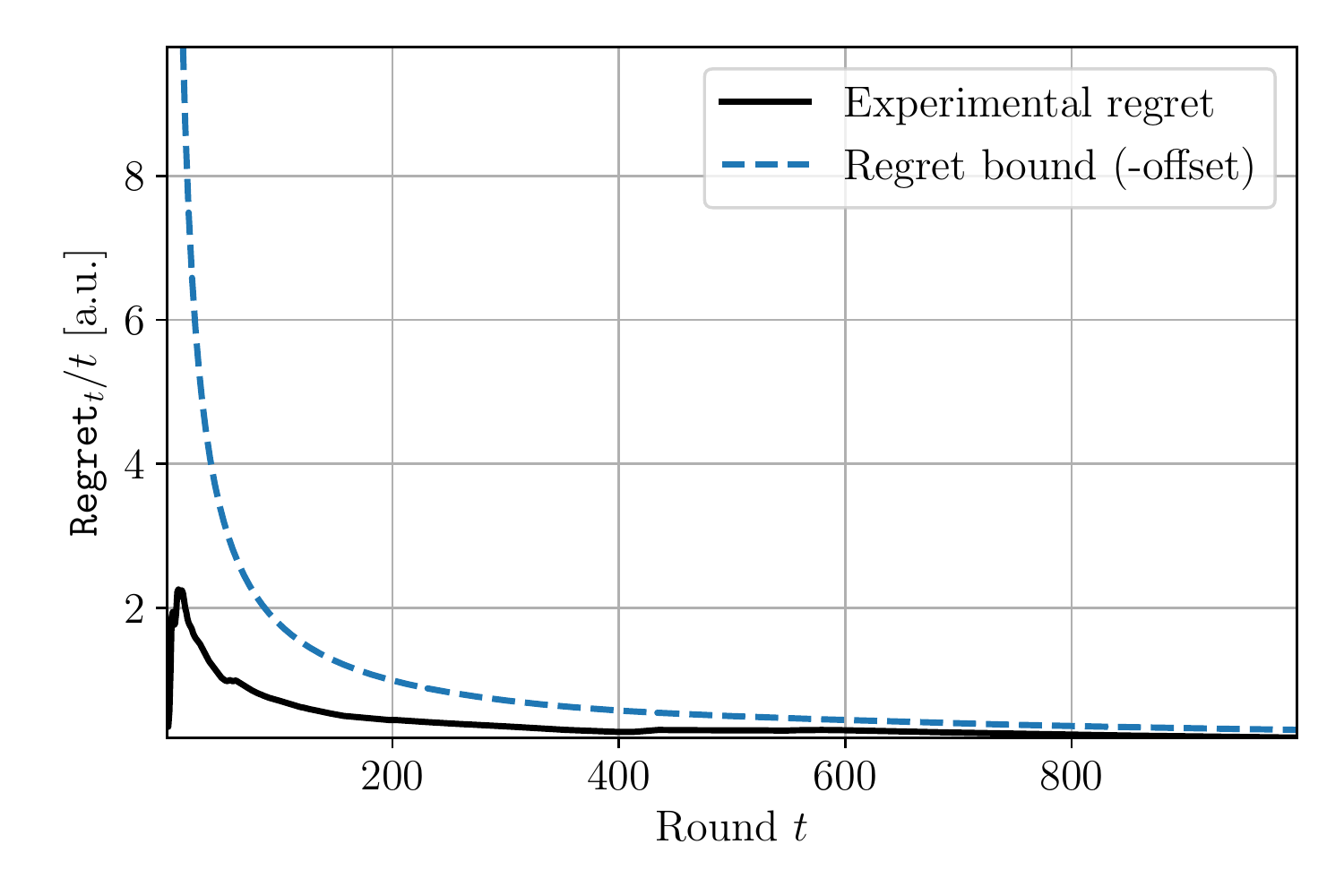}
  \vspace{-0.5cm}
  \caption{Average absolute network regret for \texttt{D-ODWDA}}
  \label{fig:avg_regret}
\end{figure}

We present next the average absolute regret at round $t$ expressed as 
$
\frac{1}{t}\sum_{\tau=1}^t \left| \sum_{i=1}^n  f_{i,\tau}\left(\xx_{i,\tau}\right) - f_\tau \left( \xx^\star_\tau \right)\right|.
$
We remark that the absolute round error upper bounds the round error. An average absolute regret that vanishes with time implies that the average regret behaves similarly and thus that the regret is sublinear. The experimental average absolute regret and the average bound are presented in Figure~\ref{fig:avg_regret}. Figure~\ref{fig:avg_regret} shows that the average absolute regret goes to zero as time increases and outperforms the bound. 

\section{Conclusion} 
In this work, we provide a dynamic regret bound for the distributed, static OCO algorithm proposed in~\cite{hosseini2016online}. \texttt{D-ODWDA} has a tighter regret bound with respect to time in comparison to all previously proposed distributed OCO algorithms. 
We consider fast timescale DR for commercial buildings with HVAC systems' air handler fan and equipped with variable frequency drive. We use \texttt{D-ODWDA} and formulate a performance-guaranteed distributed and dynamic online approach for DR of commercial buildings. The approach is scalable to large aggregations of buildings, does not require exhaustive communication infrastructure, promotes privacy, and minimizes unreliability and security risks. Lastly, we show in numerical simulations the performance of our approach to track frequency regulation signals. 

\appendix

\subsection{Proof of Lemma~\ref{lem:dual_avg_err}}
\label{app:error_y}
Using~\cite[Lemma 1]{hosseini2016online}, we have 
$\left\| \oyit - \yit \right\|_{\ast} \leq L \sum_{k=1}^{t-2} \sum_{j=1}^n \left| \left(\Pij\right)^{k} - \ppi_j \right| + 2 L$ for time-invariant $\Pij$. By \cite[Theorem 1]{anthonisse1977exponential} there exist $\nu \in \NN$ where $\nu \geq 1$ and $\gamma \in [0,1[$ such that $\left| \left(\Pij\right)^{k} - \ppi_j \right| \leq \gamma^{\left\lfloor\frac{k}{\nu}\right\rfloor}$, for all $k\geq 1$, where $\lfloor \cdot \rfloor$ is the floor operator. This leads to
\[\left\| \oyit - \yit \right\|_{\ast} \leq L \sum_{k=1}^{t-1} \sum_{j=1}^n \gamma^{\left\lceil \frac{k+1}{\nu} \right\rceil - 1} + 2L,\]
where we used the identity $\left\lfloor \frac{k}{\nu} \right\rfloor= \left\lceil \frac{k+1}{\nu} \right\rceil - 1$ which holds for all $\nu > 0$. Hence, $\left\| \oyit - \yit \right\|_{\ast}\leq \frac{L}{\gamma} \sum_{j=1}^n\sum_{k=0}^{+\infty} \gamma^{\frac{k}{\nu}}  + 2L,$ also holds. Because $\gamma < 1$, then $\gamma^{\frac{1}{\nu}} < 1$. Evaluating the geometric series, we get $\left\| \oyit - \yit \right\|_{\ast} \leq\frac{nL}{\gamma\left(1-\gamma^{\frac{1}{\nu}}\right)}  + 2L$. \hfill \IEEEQEDhere

\subsection{Proof of Lemma~\ref{lem:norm_yit}}
\label{app:yit}
The variable $\yit$ can be written as:
\begin{align}
\yit &= \sum_{j=1}^n \left(\Pij\right)^{t-1} \yy_{j,0} + \sum_{k=1}^{t-2} \sum_{j=1}^n \left(\Pij\right)^{t-1-k} \nabla f_{i,k} \left( \xx_{i,k} \right)\nonumber\\
&\qquad + \nabla f_{i,t-1} \left( \xx_{i,t-1} \right) \label{eq:y_from_0}
\end{align}
We substitute~\eqref{eq:y_from_0} in $\oyit$ and recall that $\yy_{i,0} = \mathbf{0}$ for all $i = 1,2,\ldots, n$. Because $\pi_i < 1$ for all $i$ and $\ogt = \sum_{i=1}^n \pi_i \nabla f_{i,t} \left( \xit \right)$, we can write
\[
\oyit \leq \sum_{i=1}^n  \sum_{k=1}^{t-2} \sum_{j=1}^n \left(\Pij\right)^{t-1-k} \nabla f_{i,k} \left( \xx_{i,k} \right) + \ogt.
\]
Taking the dual norm and using the triangle inequality yields
\begin{align*}
\left\|\oyit\right\|_\ast &\leq  \sum_{i=1}^n  \sum_{k=1}^{t-2} \sum_{j=1}^n  \left(\Pij\right)^{t-1-k} \left\| \nabla f_{i,k} \left( \xx_{i,k} \right) \right\|_\ast + \left\| \ogt \right\|_\ast\\
\phantom{\left\|\oyit\right\|_\ast}&\leq L \sum_{i=1}^n  \sum_{k=1}^{t-2} \sum_{j=1}^n  \left(\Pij\right)^{k}  + L,
\end{align*}
where we used $\left\| \nabla f_{i,k} \left( \xx_{i,k} \right) \right\|_\ast\leq L$ and $\left\| \ogt \right\|_\ast\leq L$ to obtain the last inequality. Let $p=\max_{i,j} \Pij$. By our assumption on the network, we have $p \in ]0,1[$. We upper bound the geometries series and obtain
$\left\|\oyit\right\|_\ast \leq \frac{n^2 L}{1-p}  + L$. \hfill \IEEEQEDhere


\end{document}